\documentclass[11pt,letterpaper]{amsart}
\usepackage{iftex}
\ifPDFTeX
  \usepackage[utf8]{inputenc}
\fi
\ifLuaTeX
  \usepackage{luatex85}
\fi
\usepackage[T1]{fontenc}
\usepackage{amsmath,amscd,amssymb}
\usepackage[mathscr]{eucal}
\usepackage{mathrsfs}
\usepackage{color,textcomp}
\usepackage{cite}
\usepackage[all,cmtip]{xy}
\IfFileExists{microtype.sty}{\usepackage{microtype}}{}
\usepackage[hypertexnames=false,hidelinks]{hyperref}
\numberwithin{equation}{section}

\newcommand{\CC}{\mathbb{C}}

\newcommand{\RR}{\mathbb{R}}
\newcommand{\ZZ}{\mathbb{Z}}

\newcommand{\cal}{\mathcal}

\def\cA{{\cal A}}

\def\cE{{\cal E}}
\def\cF{{\cal F}}

\def\cL{{\cal L}}

\def\cO{{\cal O}}

\def\cX{{\cal X}}

\DeclareMathOperator{\BC}{BC}

\DeclareMathOperator{\vol}{vol}

\newtheorem{prop}{Proposition}[section]

\newtheorem{theo}[prop]{Theorem}
\newtheorem{lemm}[prop]{Lemma}
\newtheorem{coro}[prop]{Corollary}
\newtheorem{rema}[prop]{Remark}

\newtheorem{defi}[prop]{Definition}
\newtheorem{conj}[prop]{Conjecture}

\newtheorem{defi-prop}[prop]{Definition-Proposition}

\def\dbar{\overline{\partial}}

\def\beq{\begin{equation}}
\def\eeq{\end{equation}}

\def\vol{\text{vol}}

\begin{document}
\title{Holomorphic Limits of K\"ahler Manifolds}

\author{Mu-lin Li}
\address{School of Mathematics, Hunan University, China}
\email{mulin@hnu.edu.cn}


\date{\today}

\begin{abstract}
Let $\pi:\cX\to\Delta$ be a smooth family of compact complex manifolds over the unit disc, and assume that the fibers $\cX_t$ are K\"ahler for all $t\in\Delta^*$. If $\pi$ is K\"ahler at some point $t_0\in\Delta^*$, we construct a positive $d$-closed $(1,1)$-current $T$ on the central fiber $\cX_0$ such that $\widetilde{\mathrm{vol}}_n(\{T\})>0$. This implies that $\cX_0$ belongs to Fujiki's class $\mathcal C$ provided its upper volume is finite. We also prove that the central fiber of a smooth family whose punctured fibers are projective is Moishezon, with no additional hypothesis on the central fiber.
\end{abstract}
\maketitle

\section{Introduction}

A \emph{smooth family} of compact complex manifolds will mean a proper holomorphic submersion. Let
\[
 \pi:\cX\to\Delta
\]
be such a family over the unit disc, and write $\cX_t=\pi^{-1}(t)$. Throughout the paper, the fibers are assumed to have complex dimension $n$. Suppose that $\cX_t$ is K\"ahler for every $t\in\Delta^*:=\Delta\setminus\{0\}$. A natural question is whether the central fiber $\cX_0$ must also be K\"ahler. This problem has been studied since the early development of deformation theory in complex geometry.

In complex dimension two, the answer is affirmative by Kodaira's classification of compact complex surfaces together with Siu's theorem \cite{Si83}; see also \cite{Bu,Lm99}. In higher dimensions, however, K\"ahlerness need not be preserved under specialization. Hironaka constructed a smooth deformation \cite{Hi62} whose general fibers are K\"ahler whereas the central fiber is not. Nevertheless, Hironaka's central fiber is bimeromorphic to a K\"ahler manifold. Recall that a compact complex manifold belongs to \emph{Fujiki's class $\mathcal C$} if it is bimeromorphic to a compact K\"ahler manifold. This motivates the following conjecture.
\begin{conj}\label{Fujiki}
Let $\pi:\cX\to\Delta$ be a smooth family of compact complex manifolds. If $\cX_t$ is K\"ahler for every $t\in\Delta^*$, then the central fiber $\cX_0$ belongs to Fujiki's class $\mathcal C$.
\end{conj}

More generally, one may ask which geometric properties are preserved under holomorphic specialization. For instance, if all fibers over $\Delta^*$ are projective, Moishezon, K\"ahler, or in Fujiki's class $\mathcal C$, must the central fiber have the corresponding property? Recall that a compact complex manifold is \emph{Moishezon} if it is bimeromorphic to a projective manifold. These questions are closely related, and we briefly recall several results relevant to the present work.

Popovici \cite{Po13} proved that if every fiber $\cX_t$, $t\in\Delta^*$, is projective and the central fiber $\cX_0$ satisfies either of the following conditions:
\begin{enumerate}
\item the Hodge number $h^{0,1}(\cX_t)$ is deformation invariant near $0$;
\item $\cX_0$ admits a \emph{strongly Gauduchon metric} $\omega$, namely a Hermitian metric for which $\partial\omega^{n-1}$ is $\bar\partial$-exact,
\end{enumerate}
then $\cX_0$ is Moishezon. Barlet \cite{Ba15} obtained an analogous result when the punctured fibers are Moishezon. Rao and Tsai \cite{RT21} subsequently showed that if a family contains uncountably many Moishezon fibers, then every fiber satisfying either of the above conditions is Moishezon. In forthcoming work of Rao \cite{R}, he proves that the auxiliary conditions can be removed when the family contains uncountably many Moishezon fibers. For the general fibers are K\"ahler manifolds, the author and Liu \cite{LL24} proved that the central fiber is also K\"ahler under some mild conditions of metric and cohomology on the central fiber.

Our first main result addresses Conjecture~\ref{Fujiki} by constructing a positive current of positive volume on the central fiber.
\begin{theo}[cf. Theorem~\ref{limit}]\label{projective}
Let $\pi:\cX\to\Delta$ be a smooth family of compact complex manifolds. Assume that $\cX_t$ is K\"ahler for every $t\in\Delta^*$. If $\pi$ is K\"ahler at some point $t_0\in\Delta^*$ in the sense of Definition~\ref{defi-1}, then there exists a positive $d$-closed $(1,1)$-current $T$ on $\cX_0$ such that
\[
 \widetilde{\mathrm{vol}}_n(\{T\})>0.
\]
\end{theo}

This result reduces Conjecture~\ref{Fujiki}, under the additional hypothesis that $\pi$ is K\"ahler at one point of the punctured disc, to the following conjecture proposed by Boucksom \cite{Bou}, often referred to as the Boucksom--Demailly--P\u{a}un conjecture.
\begin{conj}\label{DP}
Let $X$ be a compact complex manifold carrying a positive $d$-closed $(1,1)$-current $T$ such that $\widetilde{\mathrm{vol}}_n(\{T\})>0$. Then $X$ belongs to Fujiki's class $\mathcal C$.
\end{conj}
Boucksom--Guedj--Lu \cite[Theorem C]{BGL} established this criterion under the additional assumption that the upper volume $\overline{\vol}_n(\omega_X)$ is finite for a Hermitian metric $\omega_X$ on $X$. The upper volume, defined in \eqref{definition}, was introduced and studied by Guedj--Lu \cite{GL}; its finiteness is independent of the choice of Hermitian metric by \cite[Theorem A]{GL}. Whether the upper volume is finite on every compact complex manifold remains an important open problem. To the best of our knowledge, no compact complex manifold with infinite upper volume is currently known.

Combining Theorem~\ref{projective} with the criterion of Boucksom--Guedj--Lu yields the following consequence.
\begin{coro}
Let $\pi:\cX\to\Delta$ be a smooth family of compact complex manifolds. Assume that $\cX_t$ is K\"ahler for every $t\in\Delta^*$. Suppose moreover that
\begin{itemize}
\item $\pi$ is K\"ahler at some point $t_0\in\Delta^*$ in the sense of Definition~\ref{defi-1};
\item the central fiber $\cX_0$ admits a Hermitian metric $\omega_{\cX_0}$ such that $\overline{\mathrm{vol}}_n(\omega_{\cX_0})<\infty$.
\end{itemize}
Then $\cX_0$ belongs to Fujiki's class $\mathcal C$.
\end{coro}

Our second main result treats the projective case without any auxiliary hypothesis on the central fiber.
\begin{theo}[cf. Theorem~\ref{limit-proj}]\label{project-1}
Let $\pi:\cX\to\Delta$ be a smooth family of compact complex manifolds. If $\cX_t$ is projective for every $t\in\Delta^*$, then the central fiber $\cX_0$ is Moishezon.
\end{theo}
Thus Conjecture~\ref{Fujiki} holds, in particular, when all punctured fibers are projective.

We briefly describe the strategy of the paper. We first review the volume theory of pseudo-effective Bott--Chern classes and compare the volume of a $d$-closed class in the sense of Boucksom--Guedj--Lu with the volume defined using the absolutely continuous parts of positive currents. We then study relative K\"ahler forms. Using the deformation theorem for K\"ahler cones of Demailly--P\u{a}un \cite[Theorem 0.9]{DP04}, we construct a class $\alpha\in H^2(\cX,\RR)$ and an open dense subset $U\subset\Delta^*$ such that $\alpha_t$ is K\"ahler for $t\in U$ and of type $(1,1)$ for every $t\ne0$. Since the class need not remain of type $(1,1)$ on the central fiber, we construct $d$-closed $(1,1)$ representatives on the punctured fibers with at most a finite-order pole in the parameter. We then normalize the K\"ahler representatives by the Calabi--Yau theorem. A uniform Gauduchon mass estimate, together with Boucksom's semicontinuity theorem for the absolutely continuous parts of positive currents \cite[Proposition 2.1]{Bou}, produces a positive $d$-closed $(1,1)$-current $T$ on $\cX_0$ with $\widetilde{\mathrm{vol}}_n(\{T\})>0$. Finally, in the projective case, we combine a persistent integral Hodge class with elliptic compactness to obtain a big line bundle on the central fiber.

\subsection*{Acknowledgements}
The author is grateful to Xiao-Lei Liu and Sheng Rao for valuable discussions, and to Jian Chen for helpful comments concerning the homotopy type of Stein manifolds used in the proof of Proposition~\ref{res1}. The author was supported by the National Natural Science Foundation of China (NSFC), Grants Nos.~12271073 and 12271412. AI plays a supporting role in this project. We used the  OpenAI GPT-6.0-pro model to revise the gaps in the previous version.  The author is fully responsible for all assertions in this paper.

\section{Preliminaries on cohomology and volumes}

We collect the cohomological and pluripotential-theoretic preliminaries used throughout the paper.

Let $X$ be an $n$-dimensional compact complex manifold endowed with a Hermitian metric $\omega_X$.

\subsection{Bott--Chern and Aeppli cohomology}
For $p,q\in\mathbb Z_{\ge0}$, let $\cE^{p,q}(X)$ denote the Fr\'echet space of smooth $(p,q)$-forms on $X$, and set $\cE(X)=\cE^{0,0}(X)$. These spaces form the double complex $(\cE^{\bullet,\bullet}(X),\partial,\dbar)$, with $d=\partial+\dbar$.
The {\em Bott--Chern cohomology} of $X$ is the bigraded algebra
\beq H^{\bullet,\bullet}_{BC}(X,\CC) := \frac{\mathrm{Ker}\partial\cap\mathrm{Ker}\dbar}{\mathrm{Im}\partial\dbar}.\nonumber\eeq
Denote by $H^{1,1}_{BC}(X,\mathbb{R})\subset H^{1,1}_{BC}(X,\mathbb{C})$ the subset of classes represented by real $(1,1)$-forms. Associated with a Hermitian metric is the Bott--Chern Laplacian
$$
\Delta_{BC}:=(\partial\dbar)(\partial\dbar)^*+(\partial\dbar)^*(\partial\dbar)+(\dbar^*\partial)(\dbar^*\partial)^*+(\dbar^*\partial)^*(\dbar^*\partial)+\dbar^*\dbar+\partial^*\partial,
$$
which is self-adjoint and elliptic. We shall use the following Hodge decomposition.

\begin{prop}[{\cite[Theorem 2.3]{A1}}]
Let $X$ be a compact complex manifold endowed
with a Hermitian metric. Then there exists an orthogonal decomposition

\begin{equation}\nonumber
\cE^{p,q}(X)=\mathrm{Ker}\Delta_{BC}\oplus\mathrm{Im}\partial\dbar\oplus(\mathrm{Im}\partial^*+\mathrm{Im}\dbar^*),
\end{equation}
and an isomorphism
\begin{equation}\nonumber
H^{\bullet,\bullet}_{BC}(X,\CC)\cong \mathrm{Ker}\Delta_{BC}.
\end{equation}
\end{prop}

The {\em Aeppli cohomology} of $X$ is the bigraded $H^{\bullet,\bullet}_{BC}(X,\CC)$-module
$$ H^{\bullet,\bullet}_{A}(X,\CC) := \frac{\mathrm{Ker}\partial\dbar}{\mathrm{Im}\partial+\mathrm{Im}\dbar} .$$
Bott--Chern and Aeppli cohomology are paired by integration as follows.

\begin{theo}[{\cite[Lemma 2.5]{Sch}}]\label{pair}
Let $X$ be an $n$-dimensional compact complex manifold. For any $(p,q)\in\mathbb{N}^2$, the wedge product induces a nondegenerate bilinear pairing
\beq
\begin{CD}
  H^{p,q}_{BC}(X,\CC)\times H^{n-p,n-q}_{A}(X,\CC)@>^\wedge>>H^{n,n}_{A}(X,\CC)@>^{\int_{X}}>>\mathbb{C},\nonumber
  \end{CD}
\eeq
\beq
\begin{CD}
  (\alpha,\beta) @>>>\int_X\alpha\wedge\beta.
\end{CD}\nonumber
\eeq

\end{theo}

\subsection{Volumes of pseudo-effective classes}

A function $\varphi:X\to\RR\cup\{-\infty\}$ is called \emph{quasi-plurisubharmonic} (quasi-psh) if it is locally the sum of a plurisubharmonic function and a smooth function.
\begin{defi}
For a smooth real $(1,1)$-form $\theta$, let $PSH_{\theta}(X)$ denote the set of $\theta$-psh functions that are not identically $-\infty$; explicitly,
$$
PSH_{\theta}(X):=\{\varphi:X\to\RR\cup\{-\infty\}\mid \varphi\text{ is quasi-psh and }dd^c\varphi\ge-\theta\text{ locally}\}.\nonumber
$$
\end{defi}

When $\theta=\omega_X>0$ is a Hermitian form, define
\begin{equation}\label{definition}
\overline{\vol}_j(\omega_X)=\sup\left\{\int_X(\omega_X+dd^c\varphi)^j\wedge\omega_X^{n-j}\;\middle|\;\varphi\in PSH_{\omega_X}(X)\cap\cE(X)\right\}.
\end{equation}
For further details on these quantities, see \cite[Section 3]{GL}.

For every $B>0$, multiplication by $B$ identifies $PSH_{\omega_X}(X)$ with $PSH_{B\omega_X}(X)$. Consequently,
\beq
\overline{\vol}_n(B\omega_X)=B^n\overline{\vol}_n(\omega_X).\nonumber
\eeq

A current $T$ on $X$ is called \emph{quasi-closed} if $dT$ is smooth. By \cite[Theorem 2.2]{BGL}, the non-pluripolar product $\langle T^n\rangle$ is well defined for every quasi-closed current $T$.

Define the Bott--Chern space by
\beq
\BC^{1,1}(X):=\frac{\cE^{1,1}(X)}{\{dd^c\cE(X)\}}.\nonumber
\eeq
This is a Hausdorff Fr\'echet space, and \cite[Section 1.1]{BGL} gives
\beq
\BC^{1,1}(X)\cong \frac{\{T\, \text{quasi-closed}\, \mbox{$(1, 1)$-{current}}\}}{\{dd^cU| U\,\text{distribution}\}}.\nonumber
\eeq
For a quasi-closed current $T$, we denote its $dd^c$-class by $\{T\}\in\BC^{1,1}(X)$. A real $(1,1)$-current $T$ is called \emph{quasi-positive} if $T\ge\gamma$ for some smooth real $(1,1)$-form $\gamma$. Thus $T$ is both quasi-positive and quasi-closed if and only if $T = \theta+dd^c\varphi$ for some smooth real $(1,1)$-form $\theta$ and some quasi-psh function $\varphi$.

Let $T=\theta+dd^c\varphi$ and $T^{\prime}=\theta+dd^c\varphi^{\prime}$ be quasi-positive and quasi-closed currents. We say that $T$ is \emph{more singular} than $T^{\prime}$ if $\varphi\le \varphi^{\prime}+O(1)$, and that $T$ and $T^{\prime}$ \emph{have equivalent singularities} if $\varphi= \varphi^{\prime}+O(1)$. The current $T$ is said to \emph{have analytic singularities} if the quasi-psh function $\varphi$ has analytic singularities.

\begin{defi}[{\cite[Definition 1.14]{BGL}}]
Let $\theta$ be a smooth real $(1,1)$-form on $X$. A class $\{\theta\}\in\BC^{1,1}(X)$ is called
\begin{itemize}
\item \emph{pseudo-effective} (abbreviated psef) if it contains a positive $(1,1)$-current $T=\theta+dd^c\varphi\ge0$;
\item \emph{big} if it contains a current $T=\theta+dd^c\varphi\ge\epsilon\omega_X$ for some $\epsilon>0$.
\end{itemize}
\end{defi}

\begin{defi}[{\cite[Definition 3.1]{BGL}}]
The \emph{lower volume} and \emph{upper volume} of a quasi-closed positive $(1,1)$-current $T$ are defined by
\beq
\underline{\mathrm{vol}}_n(T):=\inf_{S}\langle S^n\rangle,\quad \overline{\mathrm{vol}}_n(T):=\sup_S \langle S^n\rangle,\nonumber
\eeq
where $S$ ranges over all quasi-closed positive $(1,1)$-currents in the same $dd^c$-class as $T$ and with equivalent singularities.
\end{defi}

\begin{defi}[{\cite[Definition 3.8]{BGL}}]\label{big-def}
For a big $dd^c$-class $\{\theta\}\in\BC^{1,1}(X)$, its \emph{lower volume} and \emph{upper volume} are defined by
 \beq
\underline{\mathrm{vol}}_n(\{\theta\}):=\sup_{0\le T\in \{\theta\}}\underline{\mathrm{vol}}_n(T),\quad \overline{\mathrm{vol}}_n(\{\theta\}):=\sup_{0\le T\in \{\theta\}}\overline{\mathrm{vol}}_n(T),\nonumber
\eeq
where $T$ ranges over all positive currents in $\{\theta\}$.
\end{defi}
These definitions extend to pseudo-effective classes by approximation:
\begin{defi}[{\cite[Definition 3.13]{BGL}}]
If $\{\theta\}\in\BC^{1,1}(X)$ is pseudo-effective, define
\beq
\underline{\mathrm{vol}}_n(\{\theta\}):=\lim_{\epsilon\to 0}\underline{\mathrm{vol}}_n(\{\theta+\epsilon\omega_X\}),\quad \overline{\mathrm{vol}}_n(\{\theta\}):=\lim_{\epsilon\to 0}\overline{\mathrm{vol}}_n(\{\theta+\epsilon\omega_X\}).\nonumber
\eeq
\end{defi}

For $d$-closed pseudo-effective classes, Boucksom--Guedj--Lu proved that the lower and upper volumes coincide.
\begin{prop}[{\cite[Proposition 3.16]{BGL}}]\label{equal}
If $\theta$ is a smooth $d$-closed $(1,1)$-form and $\{\theta\}$ is pseudo-effective, then
\beq
\underline{\mathrm{vol}}_n(\{\theta\})=\overline{\mathrm{vol}}_n(\{\theta\}).\nonumber
\eeq
\end{prop}
Accordingly, for a $d$-closed pseudo-effective class $\{\theta\}$ we write simply $\vol_n(\{\theta\})$ for this common value. Let $\cE\subset H^{1,1}_{BC}(X,\mathbb{R})$ denote the pseudo-effective cone. Extending the volume by zero outside $\cE$, one obtains the following continuity theorem.

\begin{theo}[{\cite[Theorem 3.20]{BGL}}]\label{vanishing-vol}
The volume function $\mathrm{vol}_n:\allowbreak H^{1,1}_{BC}(X,\mathbb{R}) \to \mathbb{R}_{\ge0}$ is continuous. In
particular, it vanishes outside of the big cone.
\end{theo}

We shall also use the following criterion for membership in Fujiki's class $\mathcal C$.
\begin{theo}[{\cite[Theorem C]{BGL}}]\label{cri}
Let $X$ be a compact complex manifold endowed with a Hermitian metric $\omega_X$. Then the following conditions are equivalent:
\begin{itemize}
\item $X$ belongs to Fujiki's class $\mathcal C$.
\item $\overline{\mathrm{vol}}_n(\omega_X)<\infty$ and there exists a $d$-closed positive $(1,1)$-current $T$ with
 $\mathrm{vol}_n(\{T\})>0$.
\end{itemize}
\end{theo}

For a quasi-positive $(1,1)$-current $T$ on $X$, write its Lebesgue decomposition as $T=T_{ac}+T_{sg}$, where $T_{ac}$ denotes the absolutely continuous part. If $T\ge\gamma$ for a smooth real form $\gamma$, then $T_{ac}\ge\gamma$ almost everywhere. The absolutely continuous part $T_{ac}$ may be viewed as a $(1,1)$-form with $L^1_{\mathrm{loc}}$ coefficients, and $T_{ac}^m$ is therefore defined pointwise almost everywhere for $0\le m\le n$. Following \cite[Definition 1.3]{Bou}, define
\beq
\widetilde{\vol}_n(\{\theta\}):=\sup_{0\le T\in \{\theta\}}\int_XT_{ac}^n.\nonumber
\eeq

To compare $\vol_n$ with $\widetilde{\vol}_n$ for $d$-closed pseudo-effective classes, we use the following regularization theorem of Demailly and Boucksom.

\begin{theo}[{\cite[Proposition 3.7]{Dem92},\cite[Theorem 2.4]{Bou}}]\label{regularization}
Let $T=\theta+dd^c\varphi$ be a $d$-closed positive $(1,1)$-current on a compact Hermitian manifold $(X,\omega_X)$, where $\theta$ is a smooth real $(1,1)$-form and $\varphi$ is quasi-psh. Assume that $T\ge\gamma$ for a smooth real $(1,1)$-form $\gamma$. Then there exists a sequence of functions $\varphi_k$ with analytic singularities converging to $\varphi$ such that, with $T_k:=\theta+dd^c\varphi_k$, the following hold:
\begin{itemize}
\item $T_k\to T$ weakly;
\item $T_k\ge\gamma-\epsilon_k\omega_X$, where $\epsilon_k>0$ and $\epsilon_k\to0$;
\item $T_{k,ac}\to T_{ac}$ almost everywhere on $X$
\end{itemize}
\end{theo}

The next proposition identifies these two notions of volume in the setting needed below.
\begin{prop}\label{good}
Let $X$ be an $n$-dimensional compact complex manifold with a Hermitian metric
$\omega_X$ satisfying $\overline{\mathrm{vol}}_n(\omega_X)<\infty$. Let
$\theta$ be a $d$-closed smooth real $(1,1)$-form and assume that $\{\theta\}\in H^{1,1}_{BC}(X,\RR)$ is pseudo-effective. Then
\[
 \mathrm{vol}_n(\{\theta\})=\widetilde{\mathrm{vol}}_n(\{\theta\}).
\]
\end{prop}
\begin{proof}
We prove the two inequalities separately. The finiteness assumption ensures that all volumes below are finite; this also follows from the estimates established in the first part of the argument.

Let $T=\theta+dd^c\varphi\ge0$ be an arbitrary $d$-closed positive current in
$\{\theta\}$. Apply Theorem~\ref{regularization} and write
\[
 T_k=\theta+dd^c\varphi_k,\qquad
 T_k\ge-\varepsilon_k\omega_X,\qquad
 \varepsilon_k\downarrow0,
\]
with $T_{k,ac}\to T_{ac}$ almost everywhere. Put
$S_k=T_k+\varepsilon_k\omega_X\ge0$. Resolving the analytic singularities of
$T_k$, we obtain a modification $\rho_k:Y_k\to X$ such that
\[
 \rho_k^*T_k=\beta_k+[E_k],
\]
where $\beta_k$ is smooth and $E_k$ is an effective divisor. Hence
\[
 (\rho_k^*S_k)_{ac}=\beta_k+\varepsilon_k\rho_k^*\omega_X\ge0
\]
and, since the exceptional and polar sets have Lebesgue measure zero,
\begin{equation}\label{eq:regularized-ac-mass}
 \int_X(T_{k,ac}+\varepsilon_k\omega_X)^n
 =\int_{Y_k}(\beta_k+\varepsilon_k\rho_k^*\omega_X)^n.
\end{equation}
The invariance of the upper volume under modifications and addition of divisors, together with the definition for smooth positive forms, gives
\begin{equation}\label{eq:upper-controls-ac}
 \overline{\mathrm{vol}}_n(\{\theta+\varepsilon_k\omega_X\})
 \ge \int_X(T_{k,ac}+\varepsilon_k\omega_X)^n;
\end{equation}
see \cite[Lemma 3.5 and Definition 3.8]{BGL}. Since the forms on the right-hand
side are semipositive and converge almost everywhere to $T_{ac}$, Fatou's lemma
yields
\[
 \int_XT_{ac}^n
 \le\liminf_{k\to\infty}
 \int_X(T_{k,ac}+\varepsilon_k\omega_X)^n.
\]
By Definition 3.13 of \cite{BGL} and Proposition~\ref{equal},
\[
 \lim_{k\to\infty}
 \overline{\mathrm{vol}}_n(\{\theta+\varepsilon_k\omega_X\})
 =\mathrm{vol}_n(\{\theta\}).
\]
Combining this with \eqref{eq:upper-controls-ac}, we obtain
\[
 \int_XT_{ac}^n\le\mathrm{vol}_n(\{\theta\}).
\]
Taking the supremum over all $T\ge0$ in $\{\theta\}$ gives
\begin{equation}\label{eq:tilde-le-volume}
 \widetilde{\mathrm{vol}}_n(\{\theta\})
 \le\mathrm{vol}_n(\{\theta\}).
\end{equation}
If the right-hand side vanishes, the desired equality follows immediately.

Assume now that $\mathrm{vol}_n(\{\theta\})>0$. By
Theorem~\ref{vanishing-vol}, the class $\{\theta\}$ is big. Fix $\delta>0$.
Proposition 3.12 of \cite{BGL} gives a modification $\rho:Y\to X$ and a
decomposition
\begin{equation}\label{eq:BGL-decomposition}
 \{\rho^*\theta\}=\{\omega_Y\}+\{E\},
\end{equation}
where $\omega_Y>0$ is smooth, $E$ is an effective real divisor, and
\begin{equation}\label{eq:BGL-near-volume}
 \underline{\mathrm{vol}}_n(\omega_Y)
 >\mathrm{vol}_n(\{\theta\})-\delta.
\end{equation}
We claim that $\omega_Y$ is $d$-closed. Indeed, equality
\eqref{eq:BGL-decomposition} means that
$\rho^*\theta-\omega_Y-[E]=dd^cu$ for a distribution $u$; applying $d$ and
using $d\theta=d[E]=0$ gives $d\omega_Y=0$. Thus $\omega_Y$ is a K\"ahler
form. Proposition 3.16 of \cite{BGL} therefore gives
\begin{equation}\label{eq:omegaY-volume}
 \underline{\mathrm{vol}}_n(\omega_Y)=\int_Y\omega_Y^n.
\end{equation}

Set
\[
 R:=\rho_*(\omega_Y+[E]).
\]
Then $R$ is a $d$-closed positive $(1,1)$-current on $X$. Pushing forward
\eqref{eq:BGL-decomposition} and using $\rho_*\rho^*\theta=\theta$ for a
modification shows that $\{R\}=\{\theta\}$. Since
$(\omega_Y+[E])_{ac}=\omega_Y$, Boucksom's
push-forward formula for absolutely continuous parts
\cite[Proposition 2.2]{Bou} gives
\[
 R_{ac}=\rho_*\omega_Y\qquad\text{almost everywhere}.
\]
A modification is biholomorphic outside analytic subsets of Lebesgue measure
zero. Consequently
\begin{equation}\label{eq:push-ac-mass}
 \int_XR_{ac}^n=\int_Y\omega_Y^n.
\end{equation}
Equations \eqref{eq:BGL-near-volume}--\eqref{eq:push-ac-mass} imply
\[
 \widetilde{\mathrm{vol}}_n(\{\theta\})
 \ge\int_XR_{ac}^n
 >\mathrm{vol}_n(\{\theta\})-\delta.
\]
Letting $\delta\downarrow0$ and combining the resulting inequality with \eqref{eq:tilde-le-volume} proves the proposition.
\end{proof}

\section{K\"ahler morphisms}

We next discuss relative K\"ahler structures. Even if every fiber of a smooth proper morphism $f:\cX\to B$ is K\"ahler, the total space need not be K\"ahler when the base is K\"ahler; see \cite[Example 3.9]{Big2}. The appropriate relative notion is the following.

\begin{defi}\label{defi-1}
A smooth proper morphism $f:\cX\to B$ between complex manifolds is called a \emph{K\"ahler morphism} if there exists an open cover $\{U_i\}$ of $\cX$ and smooth real-valued functions $\varphi_i$ on $U_i$ such that $\varphi_i-\varphi_j$ is the real part of a holomorphic function on $U_i\cap U_j$, and $\sqrt{-1}\partial\dbar\varphi_i$ is positive on $T_{\cX/B}|_{U_i}$.

We say that $f$ is \emph{K\"ahler at a point} $s\in B$ if its restriction over some open neighborhood of $s$ is a K\"ahler morphism.
\end{defi}

For a K\"ahler morphism, the local forms $\sqrt{-1}\partial\dbar\varphi_i$ agree on overlaps and hence glue to a global $d$-closed real $(1,1)$-form $\omega$ on $\cX$. Its restriction $\omega|_{\cX_t}$ is a K\"ahler form on every fiber. We refer to such an $\omega$ as a relative K\"ahler form.

\subsection{Sheaves of pluriharmonic functions}
We recall the sheaf-theoretic description of Bott--Chern classes that will be used below; see \cite{Va} for further details. For a complex manifold $\cX$, let $\cE^k_{\cX}$ denote the sheaf of smooth complex-valued $k$-forms and $\cE_{\cX}^{k,\ell}$ the sheaf of smooth $(k,\ell)$-forms. Let $\cE^k_{\cX,\RR}$ denote the sheaf of smooth real-valued $k$-forms, and set
\beq
\cA_{\cX}^{\ell}:=\bigoplus_{s+t=\ell,s\ge 1, t\ge1}\cE^{s,t}_{\cX}
\eeq
and let $\cA_{\cX,\RR}^{\ell}$ be its subsheaf of real forms.

Let $PH_{\cX,\RR}$ be the sheaf of real-valued pluriharmonic functions on $\cX$. We have the exact sequence
\beq\label{real-short}
\begin{CD}
0@>>>  \RR @>>> \cO_\cX@>\text{-2Im}>>PH_{\cX,\RR}@>>>0,
\end{CD}
\eeq
 and
\begingroup
\minCDarrowwidth=1.5em
\beq\label{real}
\begin{CD}
0@>>>PH_{\cX,\RR}@>>> \cE_{\cX,\RR} @>^{\sqrt{-1}\partial\dbar}>>\cE_{\cX,\RR}^{1,1}@>^{d}>> \cA_{\cX,\RR}^{3}@>^{d}>> \cA_{\cX,\RR}^{4} @>>>\cdots
\end{CD}
\eeq
\endgroup
and \eqref{real} is a fine resolution of $PH_{\cX,\RR}$; see \cite[(4.2.3)]{Va}.

Let $\overline{\cO}_{\cX}$ denote the sheaf of antiholomorphic functions and set $PH_{\cX}:=\cO_\cX+\overline{\cO}_\cX$, the sheaf of complex-valued pluriharmonic functions. Analogously, we have
\beq
\begin{CD}
0@>>>  \CC @>>> \cO_\cX\oplus\overline{\cO}_\cX@>>>PH_{\cX}@>>>0,
\end{CD}
\eeq
and
\beq\label{comp}
\begin{CD}
0@>>>PH_{\cX}@>>> \cE_{\cX} @>^{\sqrt{-1}\partial\dbar}>>\cE_{\cX}^{1,1}@>^{d}>> \cA_{\cX}^{3} @>^{d}>> \cA_{\cX}^{4} @>>>\cdots
\end{CD}
\eeq
where \eqref{comp} is a fine resolution of $PH_{\cX}$; see \cite[Theorem 2.1]{Bigo}.

In particular, Bott--Chern cohomology admits the sheaf-theoretic description
\beq\label{BC-PH}
H^{1,1}_{BC}(\cX,\CC)\cong H^1(\cX,PH_{\cX}).
\eeq
The real and complex versions are related by the commutative diagram
\beq\label{relation-1}
\begin{CD}
0@>>>  \RR @>>> \cO_{\cX}@>\text{-2Im}>>PH_{\cX,\RR}@>>>0\\
@.@VVV@VVV@VVV\\
0@>>>  \CC @>>> \cO_{\cX}\oplus\overline{\cO}_{\cX}@>(\sqrt{-1},-\sqrt{-1})>>PH_{\cX}@>>>0.
\end{CD}
\eeq

\subsection{Relative K\"ahler forms}
Let $\mathcal K_\cX:=\cE_\cX/PH_{\cX}$ and $\mathcal K_{\cX,\RR}:=\cE_{\cX,\RR}/PH_{\cX,\RR}$. A section $\Theta\in\Gamma(\cX,\mathcal K_{\cX,\RR})$ may be represented by an open cover $\{U_i\}$ and functions $\varphi_i\in\cE_{\cX,\RR}(U_i)$ satisfying $\varphi_i-\varphi_j\in PH_{\cX,\RR}(U_i\cap U_j)$. The local forms $\sqrt{-1}\partial\dbar\varphi_i$ then glue to a $d$-closed real $(1,1)$-form $\omega$ defined by
\beq
\omega|_{U_i}=\sqrt{-1}\partial\dbar\varphi_i.
\eeq
The corresponding characteristic maps fit into the commutative diagram
\beq\nonumber
\begin{CD}
\Gamma(\cX,\mathcal{K}_{\cX,\RR})@>^{\hat{c}_1}>>  H^1(\cX,PH_{\cX,\RR}) \\
@VVV@VVV\\
\Gamma(\cX,\mathcal{K}_{\cX})@>^{c_1}>>  H^1(\cX,PH_{\cX}) .
\end{CD}
\eeq

\begin{prop}\label{global}
Let $f:\cX\to B$ be a smooth family of compact complex manifolds, where $B\subset\CC$ is a pseudoconvex domain. Then $f$ is a K\"ahler morphism if and only if there exists a class $\alpha\in H^2(\cX,\RR)$ such that
\[
 \alpha_t:=\iota_t^*\alpha
\]
is a K\"ahler class on $\cX_t$ for every $t\in B$, where $\iota_t:\cX_t\hookrightarrow\cX$ denotes the inclusion.
\end{prop}
\begin{proof}
Suppose first that $f$ is K\"ahler. By Definition~\ref{defi-1}, the local forms $\sqrt{-1}\partial\dbar\varphi_i$ glue to a global $d$-closed real $(1,1)$-form $\omega$ whose restriction to each fiber is K\"ahler. Thus $\alpha=[\omega]_{\mathrm{dR}}$ has the required property.

Conversely, assume that $\alpha\in H^2(\cX,\RR)$ restricts to a K\"ahler class on every fiber. We first show that $\alpha$ admits a global $d$-closed real $(1,1)$-representative. The short exact sequence \eqref{real-short}, together with the Leray spectral sequence, yields a commutative diagram
\beq\nonumber
\begin{CD}
 @>>>H^1(\cX,PH_{\cX,\RR})@>>>H^2(\cX,\RR) @>^{\phi}>> H^2(\cX,\cO_\cX)\\
 @.@.@VVV@V^{\cong}VV\\
 @.@.H^0(B,R^2f_*\RR)@>^{\widetilde{\phi}}>>H^0(B,R^2f_*\cO_\cX).
\end{CD}
\eeq
Since $\alpha_t$ is K\"ahler, its image in $H^2(\cX_t,\cO_{\cX_t})$ vanishes for every $t\in B$. Hence $\widetilde\phi_t(\alpha_t)=0$ for all $t$, and therefore
\[
 \phi(\alpha)=0\in H^2(\cX,\cO_\cX).
\]
It follows from \eqref{real-short} and the fine resolution \eqref{real} that $\alpha$ is represented by a $d$-closed real $(1,1)$-form $\Omega$ on $\cX$.

We next lift the class of $\Omega$ to a global section of $\mathcal K_{\cX,\RR}$. Since $f_*PH_{\cX,\RR}\simeq PH_{B,\RR}$, the exact sequence
\beq\label{relation-2}
\begin{CD}
0@>>>  \RR @>>> \cO_{B}@>\text{-2Im}>>PH_{B,\RR}@>>>0
\end{CD}
\eeq
and its associated long exact sequence apply. The domain $B$ is Stein, so Cartan's Theorem~B gives $H^k(B,\cO_B)=0$ for $k\ge1$. Moreover, $B$ has the homotopy type of a CW-complex of real dimension at most one \cite[p.~39]{Mil63}; hence $H^k(B,\RR)=0$ for $k\ge2$. Consequently,
\[
 H^k(B,PH_{B,\RR})=0\qquad(k\ge1).
\]
The Leray spectral sequence therefore gives
\[
 H^1(\cX,PH_{\cX,\RR})\simeq
 \Gamma(B,R^1f_*PH_{\cX,\RR}).
\]

From
\[
0\longrightarrow PH_{\cX,\RR}\longrightarrow\cE_{\cX,\RR}
\longrightarrow\mathcal K_{\cX,\RR}\longrightarrow0
\]
and the softness of $\cE_{\cX,\RR}$, we obtain
\begin{multline}\label{sequence-3}
0\longrightarrow f_*PH_{\cX,\RR}\longrightarrow f_*\cE_{\cX,\RR}
\longrightarrow f_*\mathcal K_{\cX,\RR}\\
\longrightarrow R^1f_*PH_{\cX,\RR}\longrightarrow0.
\end{multline}
Let $\mathcal M$ denote the image of $f_*\cE_{\cX,\RR}$ in $f_*\mathcal K_{\cX,\RR}$. Thus
\[
0\longrightarrow f_*PH_{\cX,\RR}\longrightarrow f_*\cE_{\cX,\RR}
\longrightarrow\mathcal M\longrightarrow0.
\]
Since $f_*\cE_{\cX,\RR}$ is soft and $H^k(B,PH_{B,\RR})=0$ for $k\ge1$, the associated long exact sequence gives $H^k(B,\mathcal M)=0$ for $k\ge1$. Hence the exact sequence
\[
0\longrightarrow\mathcal M\longrightarrow f_*\mathcal K_{\cX,\RR}
\longrightarrow R^1f_*PH_{\cX,\RR}\longrightarrow0
\]
yields
\beq\label{sequence-4}
\begin{CD}
\Gamma(B,f_*\mathcal K_{\cX,\RR})@>>>
\Gamma(B,R^1f_*PH_{\cX,\RR})@>>>H^1(B,\mathcal M)=0\\
@V^{\cong}VV@V^{\cong}VV\\
\Gamma(\cX,\mathcal K_{\cX,\RR})@>>>H^1(\cX,PH_{\cX,\RR}).
\end{CD}
\eeq
Therefore there exists $\Theta\in\Gamma(\cX,\mathcal K_{\cX,\RR})$ whose associated Bott--Chern class is $[\Omega]$.

It remains to arrange positivity along the fibers. Fix $t_0\in B$. Since $\alpha_{t_0}$ is K\"ahler, there exists a smooth real-valued function $\psi_{t_0}$ on $\cX_{t_0}$ such that
\[
 \sqrt{-1}\bigl(\partial\dbar\Theta|_{\cX_{t_0}}
 +\partial\dbar\psi_{t_0}\bigr)
\]
is a K\"ahler form on $\cX_{t_0}$. Choose a sufficiently small ball $V_{t_0}\Subset B$ centered at $t_0$. By smooth local triviality, $\psi_{t_0}$ extends to a smooth function $\widetilde\psi_{t_0}$ on $f^{-1}(V_{t_0})$. After shrinking $V_{t_0}$, positivity persists, so
\[
 \sqrt{-1}\bigl(\partial\dbar\Theta+
 \partial\dbar\widetilde\psi_{t_0}\bigr)|_{\cX_s}
\]
is K\"ahler for every $s\in V_{t_0}$.

Choose a countable locally finite refinement $\{U_i\}$ of the resulting cover of $B$, let $\widetilde\psi_i$ be the corresponding local functions, and let $\{\rho_i\}$ be a partition of unity subordinate to $\{U_i\}$. Then
\[
 \omega:=\sqrt{-1}\partial\dbar
 \left(\Theta+\sum_i f^*\rho_i\,\widetilde\psi_i\right)
\]
is a global relative K\"ahler form. Thus $f$ is a K\"ahler morphism.
\end{proof}

\begin{rema}We note that a related result has recently been obtained by Chen \cite{Chen} via entirely different methods.
\end{rema}

\subsection{Relative cycle spaces and openness of the K\"ahler locus}
We briefly recall the facts about relative cycle spaces and variations of Hodge structure that enter the deformation argument.

Let $f:\cX\to S$ be a smooth family of $n$-dimensional compact K\"ahler manifolds. For $0\le k\le n$, denote by $C^k(\cX/S)$ the relative Barlet space parameterizing $k$-dimensional cycles contained in the fibers of $f$, and set
\[
 C(\cX/S):=\bigcup_{0\le k\le n}C^k(\cX/S).
\]
There is a natural projection $\mu:C(\cX/S)\to S$. We shall use the following properness statement.

\begin{prop}[{\cite[Proposition 3.9]{LL24}}]\label{prop-1}
Let $f:\cX\to S$ be a smooth proper morphism whose fibers are K\"ahler, and let $A$ be a connected component of $C(\cX/S)$. Then the restriction $\mu|_A:A\to S$ is proper.
\end{prop}

We also recall the numerical characterization of the K\"ahler cone due to Demailly--P\u{a}un.
\begin{theo}[{\cite[Main Theorem 0.1]{DP04},\cite[Corollary 1.3]{CT15}}]\label{Kahler}
Let $X$ be a compact K\"ahler manifold. Then its K\"ahler cone $\mathcal K(X)$ is a connected component of the set $\mathcal P$ of real $(1,1)$-classes $\alpha$ that are numerically positive on analytic cycles, namely
\[
 \int_Z\alpha^p>0
\]
for every irreducible analytic subset $Z\subset X$ of dimension $p$.
\end{theo}

We finally fix notation for the Gauss--Manin local system; see \cite[Chapter 10]{V}, \cite[Section II.1]{G68}, and \cite[Chapter 4]{CMP}. Since $\pi:\cX\to\Delta$ is a smooth proper submersion, Ehresmann's theorem makes it a smooth fiber bundle. As the disc is contractible, the bundle is smoothly trivial. Thus, after fixing $t_0\in\Delta$, there exists a diffeomorphism
\beq\label{eq:triv}
T:\cX\longrightarrow\cX_{t_0}\times\Delta
\eeq
such that $p_2\circ T=\pi$ and $T|_{\cX_{t_0}}$ is the standard inclusion into $\cX_{t_0}\times\{t_0\}$. Write $T_0:=p_1\circ T$ and set
\[
 \kappa_t:=T_0|_{\cX_t}:\cX_t\to\cX_{t_0},
 \qquad \psi_t:=\kappa_t^{-1}.
\]
The resulting maps fit into the commutative diagram
\begin{equation}\label{diagram-1}
\xymatrix{
\cX_t \ar[r]^{\iota_t} \ar@<.5ex>[rd]^{\kappa_t}  &  \cX\ar[d]^{T_0}\ar[r]^(0.35){T}\ar[rd]^(0.35){\pi}|\hole  & \cX_{t_0}\times\Delta \ar[dl]^(0.6){p_1}\ar[d]^{p_2}\\
& \cX_{t_0}\ar@<.5ex>[lu]^{\psi_t}  & \Delta
}
\end{equation}
In particular,
\beq\label{iden-1}
\psi_t^*H^k(\cX_t,\ZZ)=H^k(\cX_{t_0},\ZZ),
\qquad
\psi_t^*H^k(\cX_t,\CC)=H^k(\cX_{t_0},\CC).
\eeq
Accordingly,
\beq\label{iden-2}
E^k:=\bigcup_{t\in\Delta}\psi_t^*H^k(\cX_t,\CC)
\eeq
may be viewed as a constant local system on $\Delta$. Equivalently,
$R^k\pi_*\CC\simeq\Delta\times H^k(\cX_{t_0},\CC)$. The associated Gauss--Manin connection
\[
\nabla^{k,GM}:R^k\pi_*\CC\otimes_{\CC}\cO_\Delta
\longrightarrow R^k\pi_*\CC\otimes_{\CC}\Omega_\Delta
\]
is characterized by the property that the section obtained by parallel transport of any class in $H^k(\cX_{t_0},\CC)$ is flat.

\begin{prop}\label{res1}
Let $\pi:\mathcal X\to\Delta$ be a smooth family whose fibers $\cX_t$ are
K\"ahler for all $t\in\Delta^*$. If $\pi$ is K\"ahler at a point
$t_0\in\Delta^*$, then there exist a class $\alpha\in H^2(\cX,\RR)$ and an
open dense subset $U\subset\Delta^*$, whose complement has at most countably
many points, such that
\begin{enumerate}
\item $\alpha_t:=\iota_t^*\alpha$ is of type $(1,1)$ for every $t\ne0$;
\item $\alpha_t$ is a K\"ahler class for every $t\in U$;
\item the restricted morphism $\pi_U:\pi^{-1}(U)\to U$ is K\"ahler.
\end{enumerate}
\end{prop}
\begin{proof}
Since $\pi$ is K\"ahler at $t_0$, after shrinking to a disc
$U_0\Subset\Delta^*$ around $t_0$, Proposition~\ref{global} provides a class
$\beta\in H^2(\pi^{-1}(U_0),\RR)$ whose restriction $\beta_t$ is K\"ahler for
all $t\in U_0$. Under the differentiable trivialization \eqref{eq:triv}, the
classes $\beta_t$ form a Gauss--Manin flat section over $U_0$. Let
$\beta_{t_0}\in H^2(X_{t_0},\RR)$ be its value and define
\[
 \alpha:=T^*p_1^*\beta_{t_0}\in H^2(\cX,\RR).
\]
Then $\alpha_t=\beta_t$ on $U_0$; in particular, $\alpha_t$ is K\"ahler there.

Every fiber over $\Delta^*$ is K\"ahler. Hence the fiberwise Hodge--de Rham
spectral sequence degenerates at $E_1$. Upper semicontinuity of the Hodge
numbers, together with the constancy of the Betti numbers, implies that the
Hodge numbers are locally constant on $\Delta^*$. Grauert's theorem therefore yields the Hodge bundles and the holomorphic Hodge filtration. The flat real
class $\alpha$ defines a holomorphic section of the quotient bundle
$R^2\pi_*\CC\otimes\mathcal O_{\Delta^*}/F^1$. Its zero locus is exactly the
set of points at which $\alpha_t$ is of type $(1,1)$. Since this zero locus contains the nonempty open set $U_0$, the identity theorem gives
\begin{equation}\label{eq:alpha-H11-punctured}
 \alpha_t\in H^{1,1}(X_t,\RR)\qquad(t\ne0).
\end{equation}

Apply the deformation theorem for K\"ahler cones of Demailly--P\u{a}un
\cite[Theorem 0.9]{DP04} to $\pi|_{\cX^*}:\cX^*\to\Delta^*$. There exists a countable union $\Sigma$ of proper analytic subsets of $\Delta^*$ such that the K\"ahler cones are invariant on
$W:=\Delta^*\setminus\Sigma$ under parallel transport for the
$(1,1)$-component of the Gauss--Manin connection. Since the base is a curve,
$\Sigma$ has at most countably many points. Since the complement of a countable subset of a planar domain is path connected and $U_0\cap W\ne\varnothing$,
choose $t_*\in U_0\cap W$. The section $t\mapsto\alpha_t$ is Gauss--Manin
flat and, by \eqref{eq:alpha-H11-punctured}, stays of type $(1,1)$; hence it is
parallel for the projected $(1,1)$-connection. Since $\alpha_{t_*}$ is
K\"ahler, Demailly--P\u{a}un's theorem gives
\begin{equation}\label{eq:alpha-kahler-very-general}
 \alpha_t\text{ is K\"ahler for every }t\in W.
\end{equation}

It remains to show that the K\"ahler locus of the flat class is open. Set
\[
 U:=\{t\in\Delta^*: \alpha_t\text{ is a K\"ahler class}\}.
\]
Fix $s\in U$ and choose a small disc $V\Subset\Delta^*$ centered at $s$.
Every fiber over $V$ is K\"ahler and, by \eqref{eq:alpha-H11-punctured},
$\alpha_t$ is of type $(1,1)$ throughout $V$. The first part of the proof of
Proposition~\ref{global}, which uses only this $(1,1)$ condition together with
cohomology and base change, produces a smooth real $d$-closed $(1,1)$-form
$\Omega$ on $\pi^{-1}(V)$ representing $\alpha|_{\pi^{-1}(V)}$. Choose a
K\"ahler form $\omega_s\in\alpha_s$. Since $X_s$ is K\"ahler, the
$\partial\bar\partial$-lemma gives a real smooth function $\varphi_s$ with
\[
 \omega_s=\Omega_s+i\partial_s\bar\partial_s\varphi_s.
\]
Extend $\varphi_s$ smoothly to $\pi^{-1}(V)$. Positivity is an open condition,
so after shrinking $V$ the forms
\[
 \Omega_t+i\partial_t\bar\partial_t\widetilde\varphi_s
\]
are K\"ahler on $X_t$ for all $t\in V$. Thus $U$ is open. By
\eqref{eq:alpha-kahler-very-general}, $W\subset U$, hence
$\Delta^*\setminus U\subset\Sigma$ is at most countable and $U$ is dense.
Finally, applying Proposition~\ref{global} to each connected component of the open Riemann surface $U$ shows that $\pi_U$ is a K\"ahler morphism.
\end{proof}

\section{Proof of the main theorem}

\subsection{Auxiliary lemmas}

\begin{lemm}\label{lem:pole}
Let $a$ be a smooth real-valued function on a disc, let $b$ be a smooth
complex-valued function, and let $\ell\geq1$ be an integer. If
\begin{equation}\label{eq:m-scalar}
 m(t):=a(t)-2\operatorname{Re}\bigl(t^{-\ell}b(t)\bigr)\geq0
 \qquad(0<|t|<\varepsilon),
\end{equation}
then $b(t)=O(|t|^\ell)$. In particular, $m$ is bounded near $0$.
\end{lemm}

\begin{proof}
Expand the smooth function $b$ to order $\ell-1$ in the real variables, expressed in terms of $t$ and $\bar t$:
\begin{equation}\label{eq:taylor}
 b(t)=\sum_{p+q<\ell}c_{pq}t^p\bar t^q+O(|t|^\ell).
\end{equation}
Suppose that some coefficient of total degree less than $\ell$ is nonzero,
and let $d<\ell$ be the smallest such degree. Write $t=re^{i\theta}$.
Multiplying by $r^{\ell-d}$ and letting $r\downarrow0$ gives
\begin{equation}\label{eq:trig-limit}
 r^{\ell-d}m(re^{i\theta})\longrightarrow
 H(\theta):=-2\operatorname{Re}\!\sum_{p+q=d}
 c_{pq}e^{i(p-q-\ell)\theta}.
\end{equation}
The left-hand side is nonnegative, hence $H(\theta)\geq0$ for every $\theta$.
On the other hand, every frequency satisfies
\[
 p-q-\ell\leq d-\ell<0.
\]
Thus $H$ has zero average on the circle:
\[
 \int_0^{2\pi}H(\theta)\,d\theta=0.
\]
A continuous nonnegative function with zero average must vanish identically,
so $H\equiv0$.

The complex trigonometric polynomial in \eqref{eq:trig-limit} contains only
strictly negative frequencies, and for fixed $p+q=d$ these frequencies are
pairwise distinct. If its real part vanishes identically, all those Fourier
coefficients must vanish. This contradicts the choice of $d$.
Consequently, every coefficient with $p+q<\ell$ in \eqref{eq:taylor} is zero.
Therefore $b(t)=O(|t|^\ell)$, and \eqref{eq:m-scalar} gives $m(t)=O(1)$.
\end{proof}

\begin{lemm}\label{lem:finite-pole-representative}
Let $\pi:\cX\to\Delta$ be a smooth proper holomorphic map and assume that
$X_t$ is K\"ahler for every $t\ne0$. Let $\alpha\in H^2(\cX,\RR)$ be a class
such that $\alpha_t$ is of type $(1,1)$ for every $t\ne0$. After shrinking
$\Delta$, there exist a smooth real $d$-closed two-form $\Omega$ representing
$\alpha$, an integer $\ell\ge1$, and a smooth $(0,1)$-form $v$ on $\cX$ such
that
\begin{equation}\label{eq:finite-pole-dbar}
 \bar\partial v=t^\ell\Omega^{0,2}.
\end{equation}
Consequently, for $t\ne0$ the forms
\begin{equation}\label{eq:finite-pole-gamma}
 u_t=t^{-\ell}v_t,\qquad
 \gamma_t:=\Omega_t-d_t(u_t+\overline{u_t})
\end{equation}
are real $d_t$-closed $(1,1)$-forms satisfying
$[\gamma_t]_{\mathrm{dR}}=\alpha_t$.
\end{lemm}
\begin{proof}
Choose a smooth real closed representative $\Omega$ of $\alpha$. Since the
base is Stein and the sheaves $R^q\pi_*\mathcal O_{\cX}$ are coherent, the
Leray spectral sequence and Cartan's Theorem B give
\begin{equation}\label{eq:finite-pole-leray}
 H^2(\cX,\mathcal O_{\cX})
 \simeq\Gamma(\Delta,R^2\pi_*\mathcal O_{\cX}).
\end{equation}
Let $q$ be the section corresponding to the image of $\alpha$ under
$H^2(\cX,\RR)\to H^2(\cX,\mathcal O_{\cX})$; in the Dolbeault model this is
$[\Omega^{0,2}]_{\bar\partial}$. Over $\Delta^*$ the fibers are K\"ahler, so
$R^2\pi_*\mathcal O_{\cX}$ is locally free and cohomology and base change
identify the value $q(t)$ with the $(0,2)$ component of $\alpha_t$. Thus
$q|_{\Delta^*}=0$. Since $R^2\pi_*\mathcal O_{\cX}$ is coherent, $q$ is a
torsion section supported at the origin; after shrinking the disc there is an
integer $\ell\ge1$ such that $t^\ell q=0$.
By \eqref{eq:finite-pole-leray},
$[t^\ell\Omega^{0,2}]_{\bar\partial}=0$ on $\cX$. The Dolbeault resolution is
fine also on the noncompact total space, hence there is a smooth $(0,1)$-form
$v$ satisfying \eqref{eq:finite-pole-dbar}. Restricting to $X_t$ and using
that $t$ is constant along the fiber gives
$\bar\partial_tu_t=\Omega_t^{0,2}$ and
$\partial_t\overline{u_t}=\Omega_t^{2,0}$. Formula
\eqref{eq:finite-pole-gamma} therefore removes precisely the $(2,0)$ and
$(0,2)$ parts of $\Omega_t$. Since it differs from $\Omega_t$ by a
$d_t$-exact form, it is real, $d_t$-closed, of type $(1,1)$, and represents
$\alpha_t$.
\end{proof}

\begin{lemm}\label{lem:gauduchon-family}
Let $\pi:\cX\to\Delta$ be a smooth family of compact complex manifolds of
complex dimension $n$. Let $\Gamma_0$ be a smooth strictly positive
$(n-1,n-1)$-form on $X_0$ satisfying
$\partial_0\bar\partial_0\Gamma_0=0$. After shrinking the disc, there is a
smooth family of strictly positive forms $\Gamma_t$ on $X_t$ such that
\[
 \Gamma_t\in\mathcal A^{n-1,n-1}(X_t),\qquad
 \partial_t\bar\partial_t\Gamma_t=0,\qquad
 \Gamma_t|_{t=0}=\Gamma_0.
\]
\end{lemm}
\begin{proof}
For $n=1$ we may take $\Gamma_t=1$. Assume $n\ge2$. By the pointwise root
theorem for positive $(n-1,n-1)$-forms \cite{Michelsohn82}, there is a unique
Hermitian form $h_0$ with $h_0^{n-1}=\Gamma_0$. Extend $h_0$ to a smooth
family of Hermitian forms $h_t$ on the fibers. Gauduchon's theorem
\cite{Gauduchon77} says that for each $t$ there is, up to a positive constant,
a unique positive smooth function $f_t$ satisfying
\[
 \partial_t\bar\partial_t(f_t h_t^{n-1})=0.
\]
Normalize $f_t$ by $\int_{X_t}f_t h_t^n=1$. The equation is the kernel equation
for a smoothly varying second-order elliptic operator on functions; its
normalized kernel is one-dimensional. Standard elliptic perturbation theory
therefore gives smooth dependence of $f_t$ on $t$. Since $h_0^{n-1}$ is
already $\partial\bar\partial$-closed, $f_0$ is a positive constant. Replacing
$f_t$ by $f_t/f_0$ and setting
\[
 \Gamma_t=(f_t/f_0)h_t^{n-1}
\]
gives the required family.
\end{proof}

\begin{theo}\label{limit}
Let $\pi:\cX\to\Delta$ be a smooth family of compact complex manifolds over a
disc $\Delta$. Assume that the fibers $\cX_t$ are K\"ahler manifolds for all
$t\in\Delta^*$, and that $\pi$ is K\"ahler at a point $t_0\in\Delta^*$. Then
there exists a positive $d$-closed $(1,1)$-current $T$ on the central fiber
$\cX_0$ with
$\widetilde{\mathrm{vol}}_n(\{T\})>0$.
\end{theo}
\begin{proof}
Let $\alpha\in H^2(\cX,\RR)$ and $U\subset\Delta^*$ be supplied by
Proposition~\ref{res1}. Thus $U$ is open and dense,
$\alpha_t$ is K\"ahler for $t\in U$, and $\alpha_t$ is of type $(1,1)$ for
every $t\ne0$.

Choose a Gauduchon metric $g_0$ on $X_0$. Applying
Lemma~\ref{lem:gauduchon-family} to $\Gamma_0=g_0^{n-1}$ (with the evident
interpretation when $n=1$), and taking the positive $(n-1)$-st roots when
$n\ge2$, we obtain a smooth family of Gauduchon metrics $g_t$ such that
\begin{equation}\label{eq:family-gauduchon}
 \partial_t\bar\partial_t(g_t^{n-1})=0.
\end{equation}
Apply Lemma~\ref{lem:finite-pole-representative} to $\alpha$. We obtain
$\Omega$, $v$, and $\ell$ and, for $t\ne0$, the $d_t$-closed $(1,1)$-forms
$\gamma_t$ of \eqref{eq:finite-pole-gamma} representing $\alpha_t$.

Define smooth functions on the disc by
\begin{equation}\label{eq:main-a-b}
 a(t)=\int_{X_t}\Omega_t\wedge g_t^{n-1},\qquad
 b(t)=\int_{X_t}d_tv_t\wedge g_t^{n-1}.
\end{equation}
Then, for $t\ne0$,
\begin{equation}\label{eq:main-mass-function}
 m(t):=\int_{X_t}\gamma_t\wedge g_t^{n-1}
 =a(t)-2\operatorname{Re}\bigl(t^{-\ell}b(t)\bigr).
\end{equation}
If $t\in U$, choose any K\"ahler representative $\omega'_t\in\alpha_t$.
Since $X_t$ is K\"ahler and $\omega'_t$ and $\gamma_t$ have the same de Rham
class, the $\partial\bar\partial$-lemma gives
$\omega'_t-\gamma_t=i\partial_t\bar\partial_t\psi_t$. Hence
\begin{equation}\label{eq:main-mass-positive}
 m(t)=\int_{X_t}\omega'_t\wedge g_t^{n-1}>0
 \qquad(t\in U).
\end{equation}
The function $m$ is continuous on $\Delta^*$ and $U$ is dense, so
$m(t)\ge0$ for every $t\ne0$. Lemma~\ref{lem:pole} applied to
\eqref{eq:main-mass-function} therefore yields, after shrinking the disc,
\begin{equation}\label{eq:main-mass-bound}
 0<m(t)\le C\qquad(t\in U).
\end{equation}
Since the integral in \eqref{eq:main-mass-positive} depends only on the class $\alpha_t$, the same bound holds for every K\"ahler representative of $\alpha_t$.

We now choose representatives adapted to the Monge--Amp\`ere equation. Fix a smooth Hermitian form
$h$ on the relative tangent bundle $T_{\cX/\Delta}$ and set
\[
 \mu_t:=\frac{h_t^n}{\int_{X_t}h_t^n}.
\]
Then $\mu_t$ is a smooth family of positive volume forms with
$\int_{X_t}\mu_t=1$. The top intersection
\begin{equation}\label{eq:main-V0}
 V_0:=\int_{X_t}\alpha_t^n
\end{equation}
is independent of $t$ and is strictly positive because $\alpha_t$ is K\"ahler
on $U$. For every $t\in U$, the Calabi--Yau theorem \cite{Yau} gives a K\"ahler
form $\omega_t\in\alpha_t$ satisfying
\begin{equation}\label{eq:Yau-normalization}
 \omega_t^n=V_0\mu_t.
\end{equation}
We emphasize that no smooth dependence of $\omega_t$ on $t$ is needed. By
\eqref{eq:main-mass-bound},
\begin{equation}\label{eq:Yau-mass-bound}
 \int_{X_t}\omega_t\wedge g_t^{n-1}\le C
 \qquad(t\in U,\ |t|\ll1).
\end{equation}

Let $\iota_t:X_t\hookrightarrow\cX$ be the inclusion and put
\[
 \Theta_t:=(\iota_t)_*\omega_t.
\]
These are positive $d$-closed $(2,2)$-currents on the total space. Fix a
Hermitian metric on a relatively compact neighborhood of $X_0$ in $\cX$.
Its restrictions to the fibers are uniformly comparable with $g_t$, and
\eqref{eq:Yau-mass-bound} gives a uniform mass bound for $\Theta_t$. Choose a
sequence $t_j\in U$, $t_j\to0$. After passing to a subsequence,
\begin{equation}\label{eq:Theta-limit}
 \Theta_{t_j}\rightharpoonup\Theta
\end{equation}
for a positive $d$-closed $(2,2)$-current $\Theta$ supported on $X_0$. The
support theorem for positive closed currents supported on a smooth
hypersurface gives a unique positive $d$-closed $(1,1)$-current $T$ on $X_0$
such that
\begin{equation}\label{eq:Theta-support}
 \Theta=(\iota_0)_*T.
\end{equation}

We next identify the cohomology class of the limit. Using the differentiable trivialization
\eqref{eq:triv}, extend any closed real $(2n-2)$-form $\eta$ on $X_0$ to a
closed form $\widetilde\eta$ near the central fiber. From
\eqref{eq:Theta-limit}--\eqref{eq:Theta-support},
\[
 \int_{X_0}T\wedge\eta
 =\lim_{j\to\infty}
 \int_{X_{t_j}}\omega_{t_j}\wedge\iota_{t_j}^*\widetilde\eta
 =\langle\alpha_0,[\eta]\rangle.
\]
Thus $[T]_{\mathrm{dR}}=\alpha_0$. In particular, $\alpha_0$ admits a positive $d$-closed $(1,1)$-current representative.

It remains to show that the absolutely continuous part of $T$ has positive top self-intersection. The argument is local. Cover $X_0$ by finitely many
holomorphic submersion charts
\[
 \Phi:\mathcal V\longrightarrow D\times B,
 \qquad \pi=\operatorname{pr}_2\circ\Phi,
\]
where $D\subset\CC^n$ is a polydisc and $B\subset\Delta$ is a disc centered at $0$. After restricting to a smaller polydisc, every slice
$\mathcal V\cap X_t$ is identified holomorphically with the same domain $D$.
For a compactly supported test form $\eta$ on $D$, extend it to a test form on
$D\times B$ independent of $t$ near the origin. Then
\eqref{eq:Theta-limit} and \eqref{eq:Theta-support} imply that, under these
holomorphic identifications,
\begin{equation}\label{eq:local-current-convergence}
 \omega_{t_j}\rightharpoonup T|_D
\end{equation}
as positive closed $(1,1)$-currents on the fixed complex domain $D$.
Boucksom's semicontinuity theorem for the absolutely continuous parts of
positive currents \cite[Proposition 2.1]{Bou} applied to
\eqref{eq:local-current-convergence} gives, almost everywhere on $D$,
\begin{equation}\label{eq:Boucksom-semicont-main}
 (T_{ac})^n
 \ge \limsup_{j\to\infty}(\omega_{t_j})^n.
\end{equation}
By the normalization \eqref{eq:Yau-normalization} and the smooth dependence of
$\mu_t$,
\[
 (\omega_{t_j})^n=V_0\mu_{t_j}\longrightarrow V_0\mu_0
\]
smoothly in every such chart. Therefore \eqref{eq:Boucksom-semicont-main}
becomes
\begin{equation}\label{eq:positive-ac-density}
 (T_{ac})^n\ge V_0\mu_0\qquad\text{almost everywhere on }X_0.
\end{equation}
Integrating and using $\int_{X_0}\mu_0=1$, we obtain
\[
 \int_{X_0}(T_{ac})^n\ge V_0>0.
\]
Since $T$ itself is an admissible current in the definition of
$\widetilde{\mathrm{vol}}_n(\{T\})$, it follows that
\[
 \widetilde{\mathrm{vol}}_n(\{T\})
 \ge\int_{X_0}(T_{ac})^n
 \ge V_0>0.
\]
This completes the proof.
\end{proof}

\section{The projective case}

We begin with two consequences of the variation of Hodge structure and the exponential sequence.

\begin{lemm}\label{lem:persistent-integral}
Let $\pi:\cX\to\Delta$ be a smooth family of compact complex manifolds and
assume that $X_t$ is projective for every $t\ne0$. Then there exist an integral
class $\alpha\in H^2(\cX,\ZZ)$ and a point $t_*\in\Delta^*$ such that
\begin{enumerate}
\item $\alpha_t:=\iota_t^*\alpha$ is of type $(1,1)$ for every $t\ne0$;
\item $\alpha_{t_*}$ is the first Chern class of an ample line bundle on
$X_{t_*}$.
\end{enumerate}
\end{lemm}
\begin{proof}
Using the differentiable trivialization \eqref{eq:triv}, identify the local system $R^2\pi_*\ZZ$ with the fixed countable lattice
$\Lambda=H^2(X_{t_0},\ZZ)$. For $\lambda\in\Lambda$, let
\[
 H_\lambda=\{t\in\Delta^*: \lambda_t\in H^{1,1}(X_t)\}.
\]
Since every fiber over $\Delta^*$ is K\"ahler, the fiberwise Hodge--de Rham
spectral sequence degenerates at $E_1$. Upper semicontinuity of the Hodge numbers, together with constancy of the Betti numbers, implies local constancy of all $h^{p,q}$; Grauert's theorem then shows that the Hodge filtration is by
holomorphic subbundles. The locus $H_\lambda$ is therefore the zero locus
of the holomorphic projection of the flat section $\lambda$ to
$(R^2\pi_*\CC\otimes\mathcal O)/F^1$, hence is a closed analytic subset of
$\Delta^*$.

Let $\Lambda_{\mathrm{amp}}\subset\Lambda$ consist of those $\lambda$ whose parallel transport is ample at some point of $\Delta^*$. Since every punctured
fiber is projective,
\[
 \Delta^*=\bigcup_{\lambda\in\Lambda_{\mathrm{amp}}}H_\lambda.
\]
Since $\Delta^*$ is a Baire space and $\Lambda_{\mathrm{amp}}$ is countable, at least one analytic set $H_\lambda$ has nonempty interior. Since
$\Delta^*$ is connected, the identity theorem gives
$H_\lambda=\Delta^*$. By the definition of $\Lambda_{\mathrm{amp}}$, the
same $\lambda$ is ample at some $t_*\in\Delta^*$. Finally, under the trivialization \eqref{eq:triv}, $\lambda$ extends to an integral class
$\alpha\in H^2(\cX,\ZZ)$.
\end{proof}

\begin{lemm}\label{lem:linebundle-punctured}
Let $\pi:\cX\to\Delta$ be as in Lemma~\ref{lem:persistent-integral}. If
$\alpha\in H^2(\cX,\ZZ)$ satisfies
$\alpha_t\in H^{1,1}(X_t)$ for every $t\ne0$, then there is a holomorphic line
bundle $\cL^*$ on $\cX^*:=\pi^{-1}(\Delta^*)$ such that
\[
 c_1(\cL^*)=\alpha|_{\cX^*}.
\]
\end{lemm}
\begin{proof}
Consider the exponential sequence on $\cX^*$,
\[
 H^1(\cX^*,\mathcal O_{\cX^*}^*)\longrightarrow
 H^2(\cX^*,\ZZ)\longrightarrow H^2(\cX^*,\mathcal O_{\cX^*}).
\]
Since $\Delta^*$ is Stein, Cartan's Theorem~B and the Leray spectral sequence give
\[
 H^2(\cX^*,\mathcal O_{\cX^*})
 \simeq\Gamma(\Delta^*,R^2\pi_*\mathcal O_{\cX^*}).
\]
Since the fibers are K\"ahler on $\Delta^*$, cohomology and base change
identify the value of the image of $\alpha$ at $t$ with the image of
$\alpha_t$ in $H^2(X_t,\mathcal O_{X_t})$. This vanishes because
$\alpha_t$ is of type $(1,1)$. Hence the image of $\alpha$ in $H^2(\cX^*,\mathcal O_{\cX^*})$ vanishes. Exactness of the exponential sequence then yields the desired holomorphic line bundle.
\end{proof}

\begin{theo}\label{limit-proj}
Let $\pi:\cX\to\Delta$ be a smooth family of compact complex manifolds over
the disc $\Delta$. Assume that the fibers $\cX_t$ are projective
for all $t\in\Delta^*$. Then $\cX_0$ is Moishezon.
\end{theo}
\begin{proof}
For $n=1$, the central fiber is a compact Riemann surface and is projective.
Assume $n\ge2$ and fix the differentiable trivialization \eqref{eq:triv}.

\smallskip
\noindent\textbf{Step 1: a line bundle on the punctured family.}
By Lemma~\ref{lem:persistent-integral}, choose
$\alpha\in H^2(\cX,\ZZ)$ and $t_*\in\Delta^*$ such that $\alpha_t$ is of type
$(1,1)$ for all $t\ne0$ and $\alpha_{t_*}$ is ample. By
Lemma~\ref{lem:linebundle-punctured}, there is a line bundle
\begin{equation}\label{eq:Lstar}
 \cL^*\longrightarrow\cX^*,\qquad
 c_1(\cL^*)=\alpha|_{\cX^*}.
\end{equation}
Write $\cL_t=\cL^*|_{X_t}$. Since
$c_1(\cL_{t_*})=\alpha_{t_*}$ is an ample class, $\cL_{t_*}$ is ample.

\smallskip
\noindent\textbf{Step 2: polynomial section growth on every punctured fiber.}
By Serre vanishing, after increasing $m_0$ if necessary,
\[
 H^q(X_{t_*},\cL_{t_*}^m)=0\qquad(q>0,\ m\ge m_0),
\]
and
\begin{equation}\label{eq:P}
 P(m):=\chi(X_{t_*},\cL_{t_*}^m)
 =\frac{V}{n!}m^n+O(m^{n-1}),\qquad
 V=\int_{X_{t_*}}\alpha_{t_*}^n>0.
\end{equation}
Fix $m\ge m_0$ and put
\[
 \cF_m=\pi^\circ_*((\cL^*)^m),\qquad
 \pi^\circ=\pi|_{\cX^*}.
\]
By Grauert's direct-image theorem, $\cF_m$ is coherent. It is torsion-free: a local
section killed by a nonzero holomorphic function on the base vanishes away
from the discrete zero set of that function and hence vanishes identically.
A torsion-free coherent sheaf on the nonsingular curve $\Delta^*$ is locally
free. By upper semicontinuity, the higher cohomology groups above continue to
vanish in a neighborhood of $t_*$; Grauert's cohomology-and-base-change theorem
then gives
\[
 \operatorname{rk}\cF_m=h^0(X_{t_*},\cL_{t_*}^m)=P(m).
\]
Hence $\cF_m$ has constant rank on the connected punctured disc. For any
$s\in\Delta^*$, the exact sequence
\[
 0\to(\cL^*)^m\xrightarrow{\,t-s\,}(\cL^*)^m
 \to(\iota_s)_*\cL_s^m\to0
\]
yields, after taking direct images, an injection
\[
 (\cF_m)_s/(t-s)(\cF_m)_s\hookrightarrow H^0(X_s,\cL_s^m).
\]
Therefore
\begin{equation}\label{eq:h0lower}
 h^0(X_t,\cL_t^m)\ge P(m)
 \qquad(t\in\Delta^*,\ m\ge m_0).
\end{equation}
Choose $m_1\ge m_0$ with $P(m_1)>0$. For each $t\ne0$, choose a nonzero section of $\cL_t^{m_1}$ and denote its divisor by $Z_t$. Then
\begin{equation}\label{eq:R_t}
 R_t:=\frac1{m_1}[Z_t]\ge0,\qquad
 \{R_t\}_{BC}=c_1^{BC}(\cL_t).
\end{equation}
No continuity of the divisors $Z_t$ is used.

\smallskip
\noindent\textbf{Step 3: finite-order-pole representatives and bounded positive pairings.}
Apply Lemma~\ref{lem:finite-pole-representative} to the image of $\alpha$ in $H^2(\cX,\RR)$. We obtain $\Omega$, $v$, $\ell$ and real $d_t$-closed $(1,1)$-forms
\begin{equation}\label{eq:gamma-projective}
 \gamma_t=\Omega_t-d_t(t^{-\ell}v_t+\bar t^{-\ell}\overline{v_t}),
 \qquad [\gamma_t]_{\mathrm{dR}}=\alpha_t.
\end{equation}
On the projective fiber $X_t$, the $\partial\bar\partial$-lemma and
\eqref{eq:Lstar} imply
\begin{equation}\label{eq:gamma-projective-class}
 \{\gamma_t\}_{BC}=c_1^{BC}(\cL_t).
\end{equation}
Fix a smooth strictly positive $\partial_0\bar\partial_0$-closed $(n-1,n-1)$-form $\Gamma_0$ on $X_0$, and extend it to a smooth family $\Gamma_t$ by Lemma~\ref{lem:gauduchon-family}. From
\eqref{eq:R_t} and \eqref{eq:gamma-projective-class},
\begin{equation}\label{eq:test-positive}
 m_\Gamma(t):=\int_{X_t}\gamma_t\wedge\Gamma_t
 =\int_{X_t}R_t\wedge\Gamma_t\ge0
 \qquad(t\ne0).
\end{equation}
Writing
\begin{align}\label{eq:test-pole}
 m_\Gamma(t)
 &=a_\Gamma(t)-2\operatorname{Re}\bigl(t^{-\ell}b_\Gamma(t)\bigr),\\
 a_\Gamma(t)&=\int_{X_t}\Omega_t\wedge\Gamma_t,\qquad
 b_\Gamma(t)=\int_{X_t}d_tv_t\wedge\Gamma_t,\nonumber
\end{align}
The functions $a_\Gamma$ and $b_\Gamma$ are smooth at $0$. Lemma~\ref{lem:pole}
therefore gives
\begin{equation}\label{eq:test-bound}
 0\le\int_{X_t}\gamma_t\wedge\Gamma_t\le C_\Gamma
 \qquad(0<|t|\ll1).
\end{equation}

Choose a Gauduchon metric $g_0$ on $X_0$. Apply Lemma~\ref{lem:gauduchon-family} to $g_0^{n-1}$; Michelsohn's root theorem \cite{Michelsohn82} then yields a smooth family of Gauduchon metrics $g_t$. Let
$r=\dim_\RR H_A^{n-1,n-1}(X_0,\RR)$ and choose real
$\partial_0\bar\partial_0$-closed representatives
$\eta_1,\ldots,\eta_r$ of an Aeppli basis. For sufficiently small
$\varepsilon_a>0$, the forms
\[
 \Gamma_{0,0}=g_0^{n-1},\qquad
 \Gamma_{a,0}=g_0^{n-1}+\varepsilon_a\eta_a\quad(1\le a\le r)
\]
are strictly positive. Their Aeppli classes span
$H_A^{n-1,n-1}(X_0,\RR)$. Extend each of them by
Lemma~\ref{lem:gauduchon-family}; estimate \eqref{eq:test-bound} holds for this
finite family.

\smallskip
\noindent\textbf{Step 4: uniform Bott--Chern harmonic representatives.}
For $t\ne0$, let $\theta_t$ be the real Bott--Chern harmonic representative
of $c_1^{BC}(\cL_t)$ with respect to $g_t$. By
\eqref{eq:gamma-projective-class},
\begin{equation}\label{eq:finite-pairing}
 \left|\int_{X_t}\theta_t\wedge\Gamma_{a,t}\right|
 =\left|\int_{X_t}\gamma_t\wedge\Gamma_{a,t}\right|
 \le C_a\qquad(0\le a\le r).
\end{equation}
We claim that for every $k\ge0$ there exists a constant $C_k$ such that
\begin{equation}\label{eq:Ckbound}
 \|\theta_t\|_{C^k(X_t,g_t)}\le C_k
 \qquad(0<|t|\ll1).
\end{equation}
Identify the fibers differentiably by \eqref{eq:triv}. The Bott--Chern
Laplacians $\Delta_{BC,t}$ form a smooth family of fourth-order elliptic
operators and satisfy uniform estimates
\begin{equation}\label{eq:ellipticBC}
 \|\xi\|_{H^{s+4}}
 \le C_s\bigl(\|\Delta_{BC,t}\xi\|_{H^s}+\|\xi\|_{L^2}\bigr).
\end{equation}
Suppose that $\|\theta_t\|_{L^2}$ were unbounded. Choose $t_j\to0$ such that
$c_j:=\|\theta_{t_j}\|_{L^2}\to\infty$ and put
$\xi_j=c_j^{-1}\theta_{t_j}$. Uniform elliptic estimates and a diagonal argument yield, after passing to a subsequence,
\[
 \xi_j\to\xi_0\quad\text{in }C^\infty,\qquad
 \Delta_{BC,0}\xi_0=0,\qquad
 \|\xi_0\|_{L^2}=1.
\]
Dividing \eqref{eq:finite-pairing} by $c_j$ and passing to the limit gives
\[
 \int_{X_0}\xi_0\wedge\Gamma_{a,0}=0\qquad(0\le a\le r).
\]
The classes $[\Gamma_{a,0}]_A$ span $H_A^{n-1,n-1}(X_0,\RR)$ and the
Bott--Chern--Aeppli pairing is nondegenerate by Theorem~\ref{pair}. Hence
$[\xi_0]_{BC}=0$. Since $\xi_0$ is Bott--Chern harmonic, this forces
$\xi_0=0$, contradicting its unit $L^2$ norm. Thus the $L^2$ norms are
bounded, and \eqref{eq:ellipticBC} plus Sobolev embedding proves
\eqref{eq:Ckbound}.

Choose a sequence $t_j\to0$ and pass to a subsequence such that
\begin{equation}\label{eq:theta-limit}
 \theta_{t_j}\to\theta_0\qquad\text{in }C^\infty.
\end{equation}
The limit is a real $d$-closed $(1,1)$-form on $X_0$. Pairing with fixed
closed $(2n-2)$-forms under \eqref{eq:triv} shows
\begin{equation}\label{eq:theta-class}
 [\theta_0]_{\mathrm{dR}}=\alpha_0.
\end{equation}

\smallskip
\noindent\textbf{Step 5: a limiting holomorphic line bundle.}
For each $j$, choose a Hermitian metric on $\cL_{t_j}$ whose Chern connection
$D_j$ has normalized curvature
\begin{equation}\label{eq:curv-normalized}
 \frac{i}{2\pi}F_{D_j}=\theta_{t_j}.
\end{equation}
Such a metric exists because $\theta_{t_j}$ represents $c_1^{BC}(\cL_{t_j})$. Let $(E,h)$ be a fixed smooth Hermitian line bundle on
the differentiable manifold $M=X_0$ with $c_1(E)=\alpha_0$. Since complex
line bundles are classified topologically by their first Chern class, the
bundles $\cL_{t_j}$, transported by \eqref{eq:triv}, are smoothly isomorphic
to $E$; choose unitary identifications and regard $D_j$ as unitary connections
on $(E,h)$.

Fix a reference unitary connection $D_{\mathrm{ref}}$ and write
\[
 D_j=D_{\mathrm{ref}}+ia_j,\qquad a_j\in\mathcal A^1(M,\RR).
\]
The curvature identity gives uniform $C^k$ bounds for $da_j$. Fix a
Riemannian metric on $M$, let $G$ be the Green operator for its Hodge
Laplacian, and set $\zeta_j=d^*G(da_j)$. Since $da_j$ is exact and closed,
Hodge theory gives $d\zeta_j=da_j$, while elliptic estimates bound
$\zeta_j$ in every $C^k$ norm. Hence
\[
 a_j=\zeta_j+df_j+\tau_j,
 \qquad \tau_j\in\mathcal H^1(M,\RR).
\]
A unitary gauge transformation removes $df_j$. The harmonic logarithmic
derivatives of $S^1$-valued gauge transformations form the full period lattice
$2\pi H^1(M,\ZZ)$ inside $\mathcal H^1(M,\RR)$, so a further gauge places
$\tau_j$ in a fixed bounded fundamental domain. After these gauges, $a_j$ is
bounded in every $C^k$ norm. Passing to a subsequence,
\begin{equation}\label{eq:Dlimit}
 D_j\to D_0\qquad\text{in }C^\infty.
\end{equation}

Let $J_{t_j}$ denote the transported complex structures on $M$. The operators
\[
 \bar\partial_j=(D_j)^{0,1}_{J_{t_j}}
\]
define holomorphic line bundles isomorphic to $\cL_{t_j}$, and they converge
to $\bar\partial_0=(D_0)^{0,1}_{J_0}$. From
\eqref{eq:curv-normalized} and \eqref{eq:theta-limit},
$F_{D_0}=-2\pi i\theta_0$ is of type $(1,1)$, so
$\bar\partial_0^2=0$. Thus $\bar\partial_0$ defines a holomorphic line bundle
$\cL_0$ on $X_0$ with $c_1(\cL_0)=\alpha_0$.

\smallskip
\noindent\textbf{Step 6: transfer of section growth to the central fiber.}
Fix $m\ge1$. On the fixed smooth bundle $E^m$ consider
\[
 P_{j,m}=\bar\partial_{j,m}^*\bar\partial_{j,m}.
\]
Define the adjoints using smoothly varying Hermitian metrics.
The nonnegative elliptic operators $P_{j,m}$ converge in $C^\infty$ to
$P_{0,m}$ and
\[
 \ker P_{j,m}=H^0(X_{t_j},\cL_{t_j}^m),\qquad
 \ker P_{0,m}=H^0(X_0,\cL_0^m).
\]
Elliptic compactness gives
\begin{equation}\label{eq:kernel-semicont}
 h^0(X_0,\cL_0^m)
 \ge\limsup_{j\to\infty}h^0(X_{t_j},\cL_{t_j}^m).
\end{equation}
Indeed, if the right-hand side is at least $N$, choose $N$ orthonormal kernel sections. Uniform elliptic estimates and Rellich compactness give, after a
subsequence, $N$ orthonormal limits in $\ker P_{0,m}$.

Combining \eqref{eq:h0lower}, \eqref{eq:P}, and
\eqref{eq:kernel-semicont}, for every $m\ge m_0$ we have
\[
 h^0(X_0,\cL_0^m)
 \ge P(m)=\frac{V}{n!}m^n+O(m^{n-1}).
\]
Thus $\cL_0$ is big; equivalently, $\kappa(X_0,\cL_0)=n$. The classical Moishezon criterion therefore implies that $X_0$ is Moishezon.
\end{proof}

\bibliographystyle{amsplain}

\end{document}